\documentclass[12pt]{amsart}
\usepackage{epsfig,color}

\makeatother

\theoremstyle{definition}

\def\fnum{equation}

\numberwithin{equation}{section}
\newcommand{\Vol}{{\text{Vol}}}
\newcommand{\V}{{\text{V}}}

\newcommand{\nn}{{\bf{n}}}

\def\RR{{\bold R}}

\def\SS{{\bold S}}

\def\CC{{\bold C }}

\newcommand{\e}{{\text {e}}}

\newcommand{\Area}{{\text {Area}}}

\def\bH{{\bold H}}

\newcommand{\cF}{{\mathcal{F}}}
\newcommand{\cL}{{\mathcal{L}}}

\newcommand{\cP}{{\mathcal{P}}}

\newcommand{\eqr}[1]{(\ref{#1})}

\title[In search of stable geometric structures]{In search of stable geometric structures}

\author{Tobias Holck Colding}%
\address{MIT, Dept. of Math.\\
77 Massachusetts Avenue, Cambridge, MA 02139-4307.}
\author{William P. Minicozzi II}%
\thanks{The  authors
were partially supported by NSF Grants DMS 1812142 and DMS 1707270.}

\email{colding@math.mit.edu and minicozz@math.mit.edu}

\begin{document}

\maketitle

\begin{abstract}
We will look for stable structures in four situations and discuss what is known and unknown.  
\end{abstract}

\section{Form finding}

Frei Otto was one of the most emblematic architects and engineers of the 20th century, renowned for his research in lightweight tensile structures.   In his book:  ``Finding form'' Frei Otto discussed the application of the optimal form in architecture:
\begin{quote}
``Natural structures are optimized, having maximum strength for minimum materials.'' 
\end{quote}

\begin{figure}[htbp]
    \begin{minipage}[t]{0.5\textwidth}
    \centering\includegraphics[totalheight=.32\textheight, width=.9\textwidth]{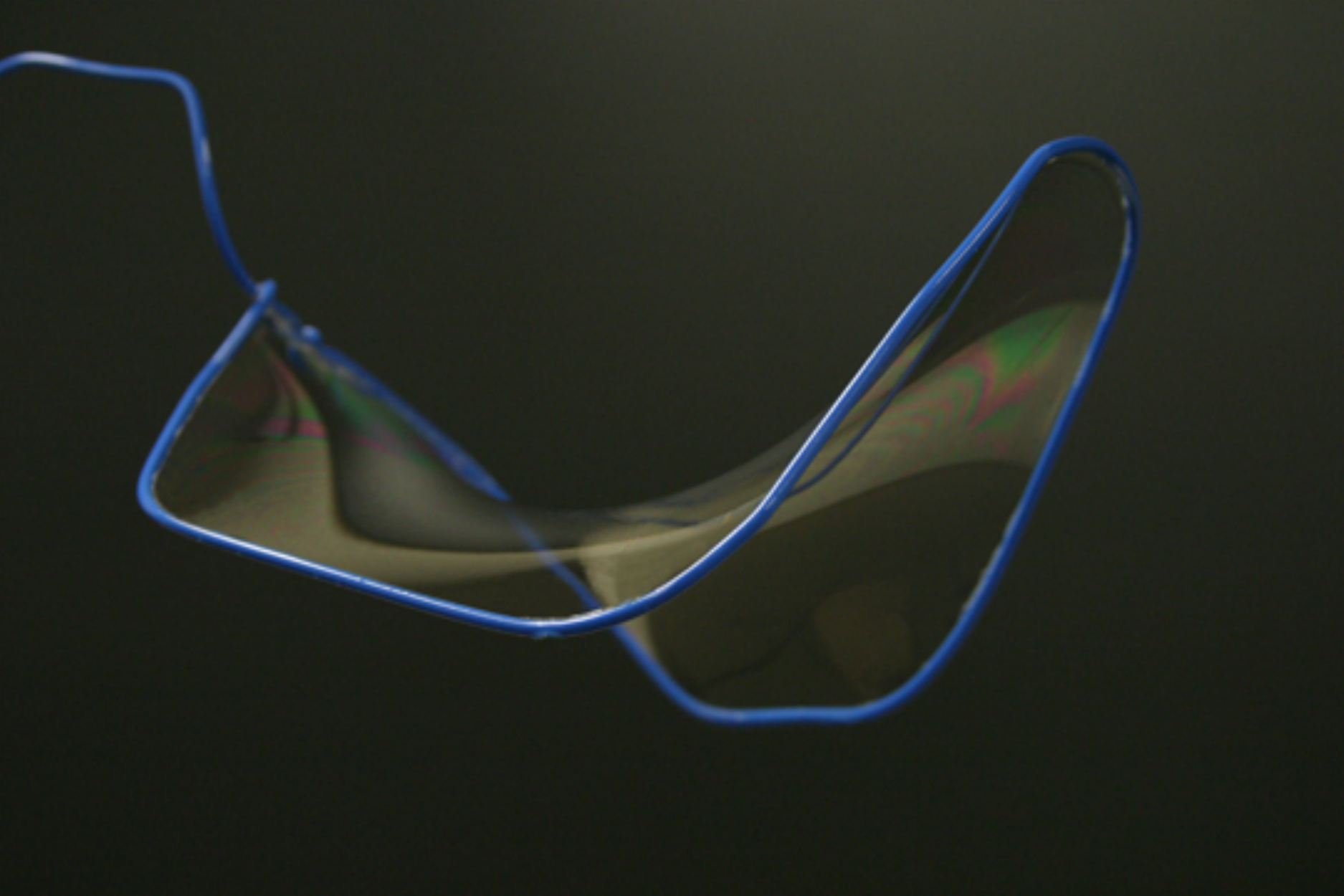}
    %\caption{The soap film minimizes surface tension.}  
    \end{minipage}\begin{minipage}[t]{0.5\textwidth}
    \centering\includegraphics[totalheight=.32\textheight, width=.9\textwidth]{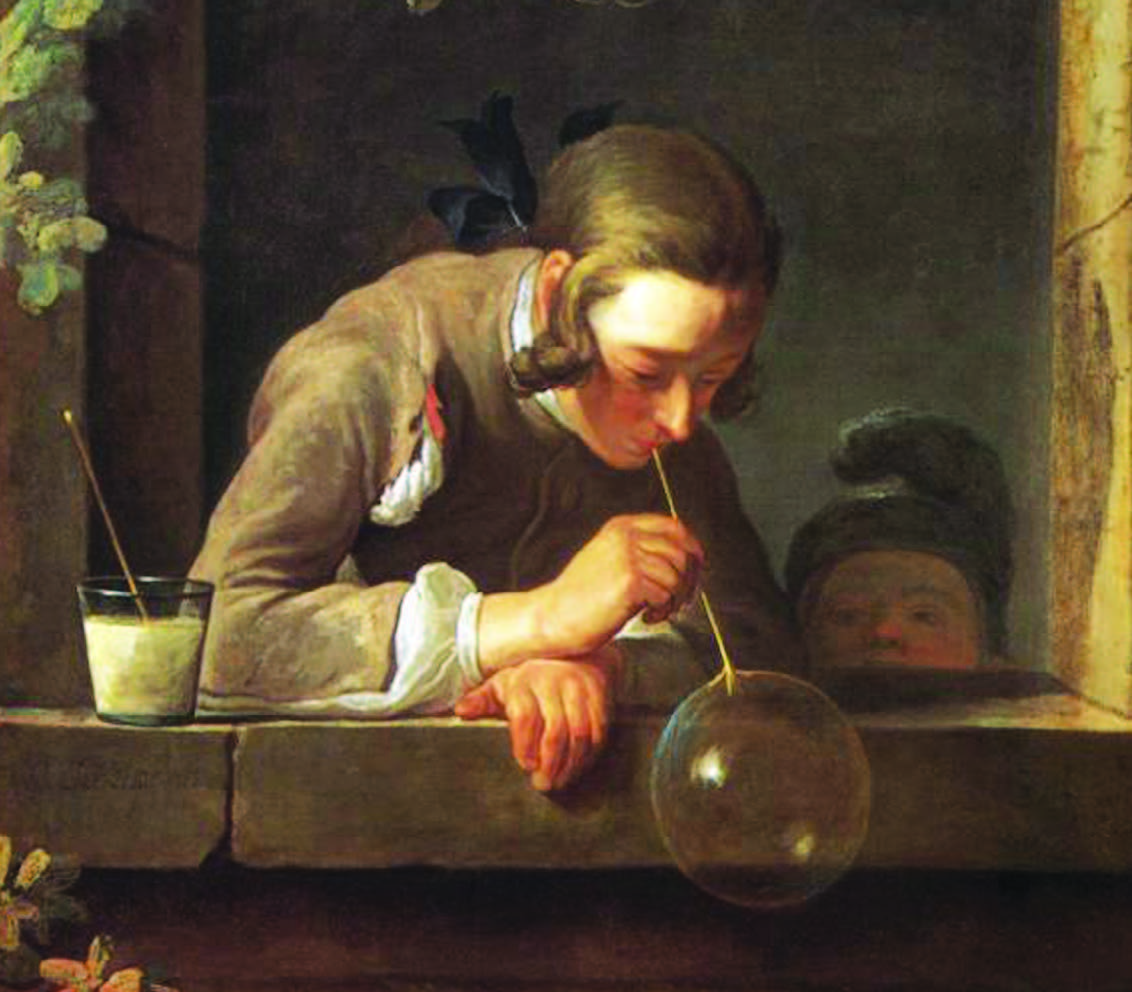}
    %\caption{A soap film on a coil.} 
    \end{minipage}
    \caption{Soap films minimize surface area for a fixed boundary.  Soap bubbles minimize area for a fixed enclosed volume.}   \label{f:pl}
\end{figure}

Natural forces tend towards stable structures like a soap film that is pulled tight by the force of surface tension. 
The soap film finds a stable equilibrium where it has the least area of any nearby surface with the same boundary.  
The  $19$-century Belgian physicist Joseph Plateau studied this experimentally, and the existence problem became known as the Plateau problem (see figure \ref{f:pl}):

%Plateau problem   (the geodesic problem one dimension up):
\begin{quote}
Given a closed curve $\gamma$, find the surface $\Gamma$   bounding $\gamma$ with the least   area.
\end{quote}
The solution $\Gamma$ is a minimal surface.  
%The problem of finding $\Gamma$ is named for the $19$-century Belgian physicist Joseph Plateau who formulated it in his study of soap films.   
The Plateau problem was finally solved around 1930 by Douglas and Rado, who worked independently.   Douglas received one of the first two Fields medals for his solution in  1936.

 \begin{figure}[htbp]
    \begin{minipage}[t]{0.5\textwidth}
    \centering\includegraphics[totalheight=.32\textheight, width=.9\textwidth]{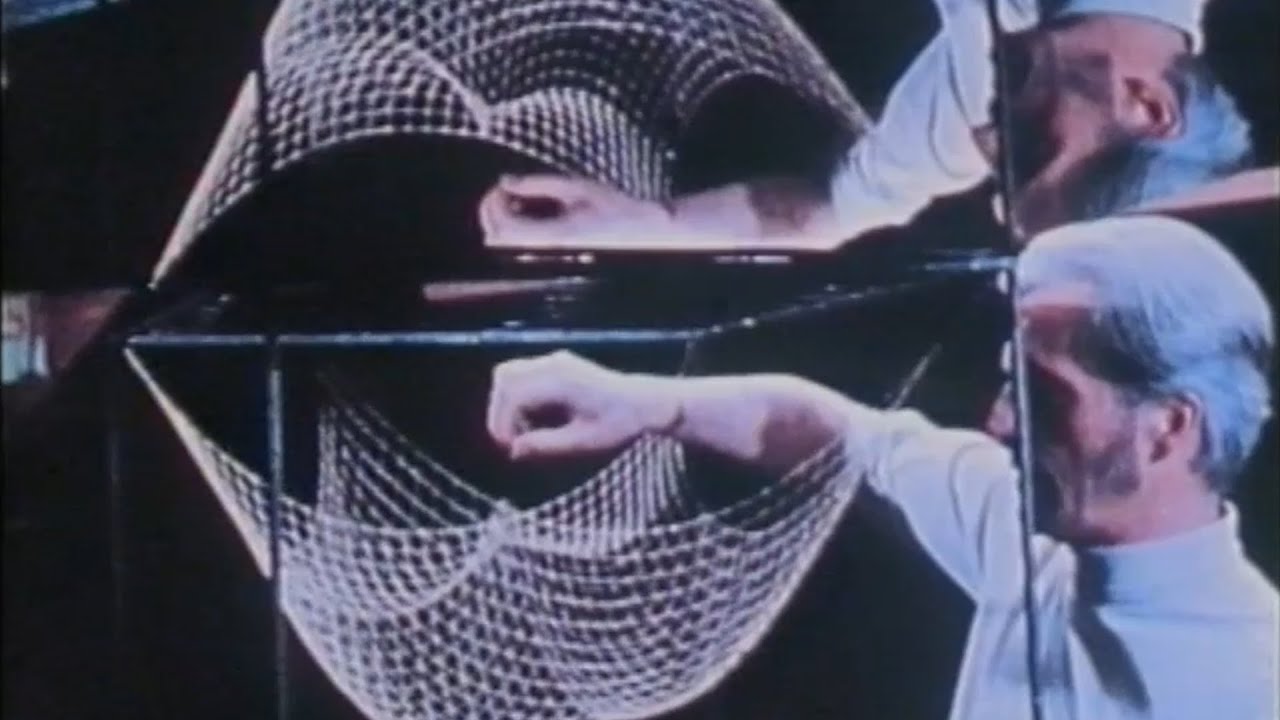}
    %\caption{The soap film minimizes surface tension.}  
    \end{minipage}\begin{minipage}[t]{0.5\textwidth}
    \centering\includegraphics[totalheight=.32\textheight, width=.9\textwidth]{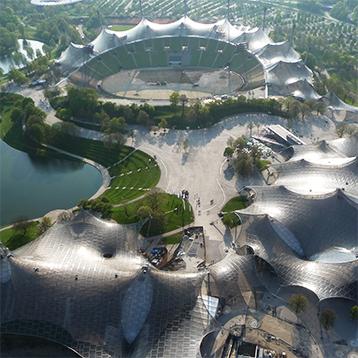}
    %\caption{A soap film on a coil.} 
    \end{minipage}
    \caption{Architect and structural engineer Frei Otto pioneered the concept of ``form finding''.  Here he is experimenting with stable structures.  To the right his famous olympic stadium based on his experiments with soap films.}   
\end{figure}

The classical Plateau problem 
searches for the minimum of the area on the space of mappings of the disk with fixed boundary.  
The space of mappings is infinite dimensional, introducing serious difficulties, but one can 
define the derivative (or first variation) of area and the notion of a critical point.  The critical points, known as minimal surfaces, were introduced by Euler and Lagrange in the 1700s.  

\section{Finding stable structures}

In mathematics, structural stability is a fundamental property of a dynamical system which means that the qualitative behavior of the trajectories is unaffected by small perturbations.   
Given a smooth function $f$ on a finite dimensional space, the gradient $\nabla f$ points in the direction of steepest ascent.
The critical points of $f$ are the points where 
  $\nabla f$ vanishes.  If $p$ is a local minimum of $f$, then
the   second derivative test  tells us that   the Hessian matrix  of $f$ at $p$ is nonnegative.  More generally, the number of negative eigenvalues of the Hessian is called the index of the critical point.   A fundamental method to find the minimum of $f$ is the method of gradient descent.  Here, we make an initial guess $p_0$ and then iteratively move in the negative gradient direction, the direction of steepest descent, by setting $p_{i+1} = p_i - \nabla f (p_i)$.  This can also be done continuously by defining a negative gradient flow
\begin{align}
	\frac{d x}{d t} = - \nabla f (x(t)) \, .
\end{align}
The function $f(x(t))$ decreases as efficiently as possible as $x(t)$  heads towards the minimum.  The dynamics near a non-degenerate critical point are determined by the index.  If the index is zero, then the critical point is attracting   and 
the entire neighborhood flows towards the critical point.  However, when the index is positive,  
a generic point   will flow out of the neighborhood, missing the critical point.

We will look for stable structures in four situations and discuss what is known and unknown.  Those four are:   Minimal hypersurfaces and minimal submanifolds of higher codimension.  Singularities that are stable or generic, and cannot be perturbed away, for motion by mean curvature of hypersurfaces and, finally, such singularities for motion by mean curvature in higher codimension.

      \begin{figure}[htbp]
\centering\includegraphics[totalheight=.18\textheight, width=1\textwidth]{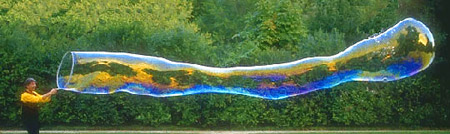}  
\caption{ Surface tension pulls the bubble towards the shape that gives the minimum surface energy - the lowest ratio of surface area to volume.}    \label{f:shrinker9a}
  \end{figure}

\section{Minimal surfaces}

 \begin{figure}[htbp]
    \begin{minipage}[t]{0.5\textwidth}
    \centering\includegraphics[totalheight=.45\textheight, width=.9\textwidth]{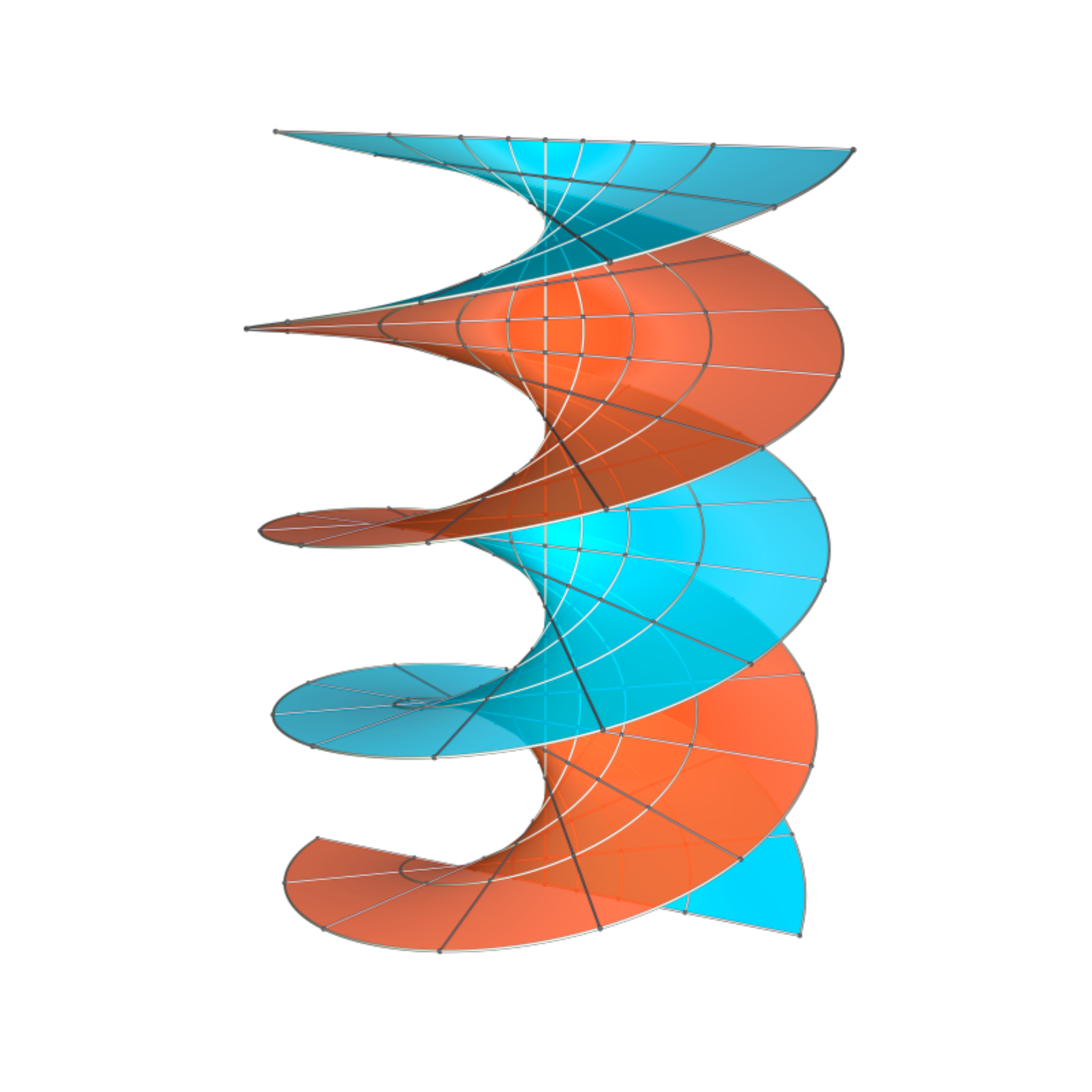}
    
    \end{minipage}\begin{minipage}[t]{0.5\textwidth}
    \centering\includegraphics[totalheight=.4\textheight, width=.9\textwidth]{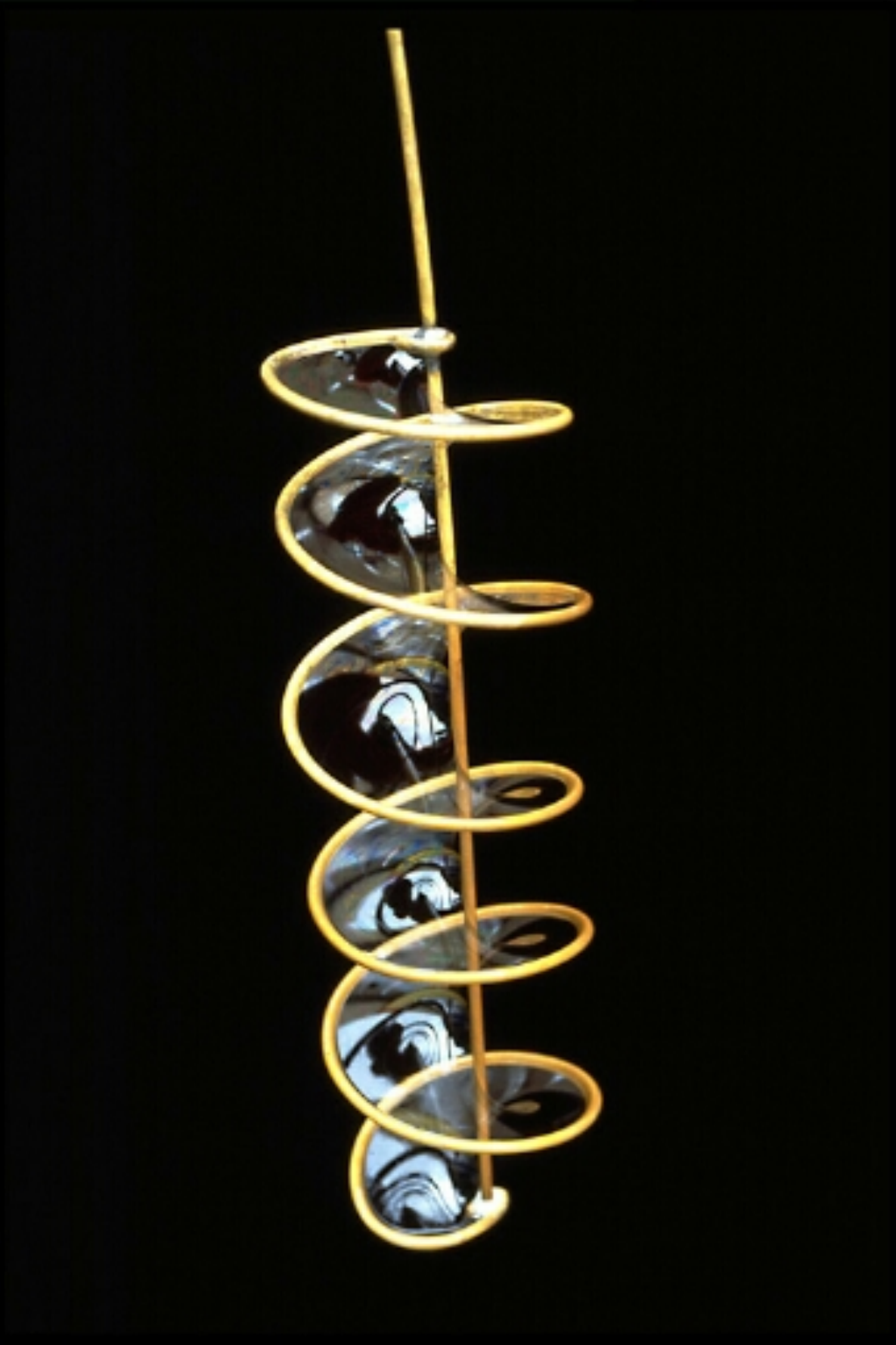}
    \end{minipage}   
    \caption{Helicoid: A minimal  surface  discovered by Meusnier in 1776.  The helicoid is a ruled surface.  Half of a helicoid is stable and can be obtained as a soap film.}  
\end{figure}

 \begin{figure}[htbp]
    \begin{minipage}[t]{0.5\textwidth}
    \centering\includegraphics[totalheight=.45\textheight, width=.9\textwidth]{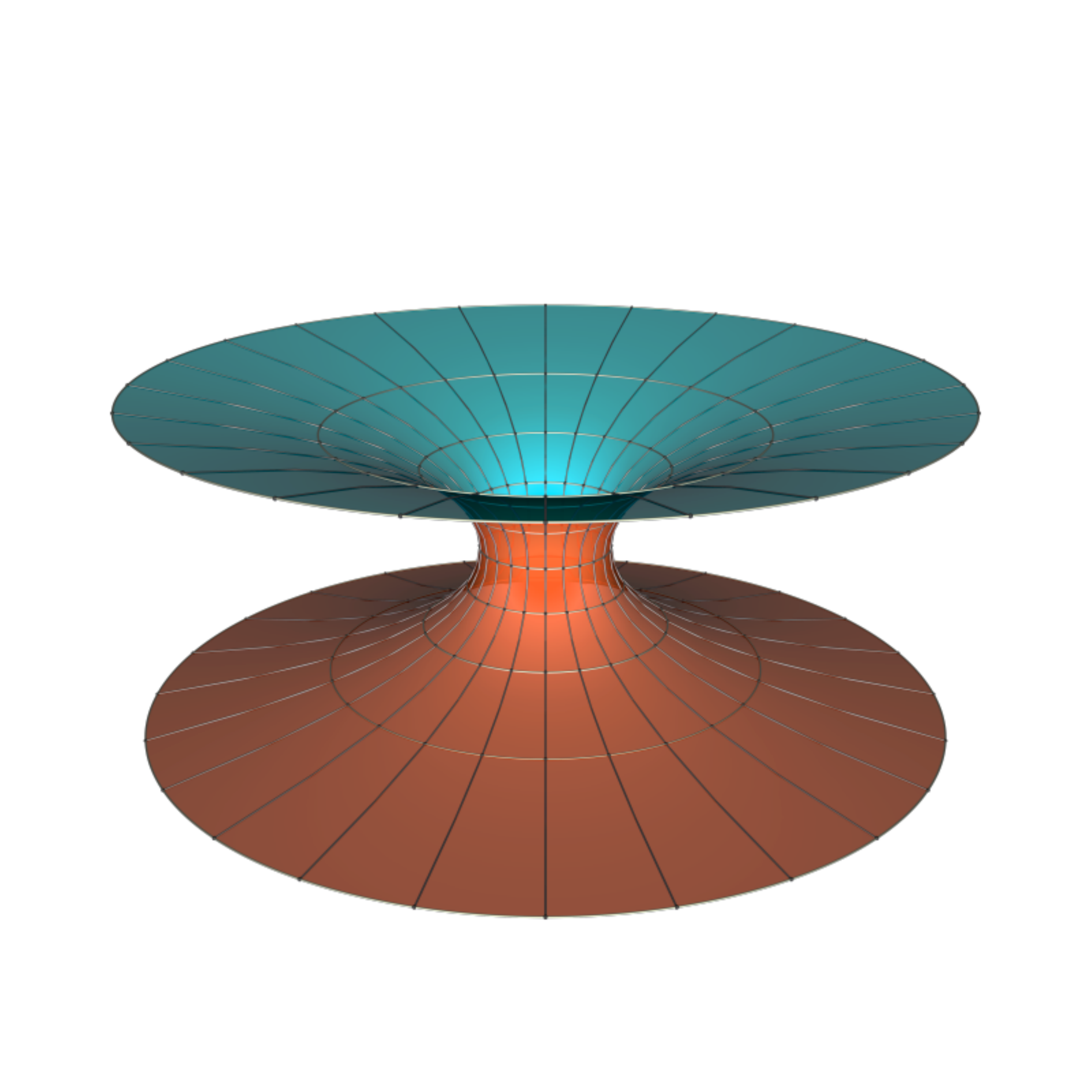}
    \end{minipage}\begin{minipage}[t]{0.5\textwidth}
    \centering\includegraphics[totalheight=.35\textheight, width=.9\textwidth]{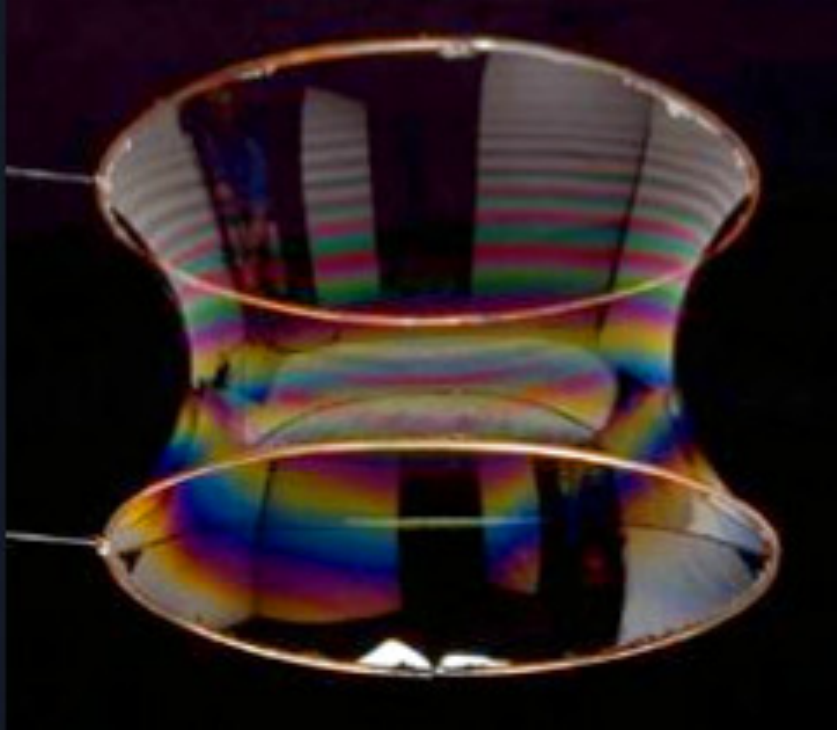}
    \end{minipage}
    \caption{The catenoid, discovered by Euler in 1741: the only minimal surface of revolution.  Part of a catenoid is stable and can be obtained as a soap film.} 
\end{figure}

Let $\Sigma^n\subset \RR^N$ be a smooth submanifold (possibly with boundary).   Given an infinitely differentiable (i.e. smooth), compactly supported, normal (orthogonal to $\Sigma$) vector field $V$ on $\Sigma$, consider the one-parameter variation
\begin{align}
\Sigma_{s,V}=\{x+s\,V(x)\,\,|\,\,x\in\Sigma\}\, .
\end{align}
This gives a path $s \to \Sigma_{s,V}$ in the space of submanifolds with $\Sigma_{0,V} = \Sigma$.
The so-called first variation formula of area or volume is the equation (integration is with respect to $d\,\Vol$)
\begin{align}   \label{e:firstvariation}
\frac{d}{ds}_{s=0}\Vol(\Sigma_{s,V})=\int_{\Sigma}\langle V,\bH\rangle\, ,
\end{align}
where  $\bH$ is the mean curvature vector.   When $\Sigma$ is a hypersurface,  $\bH$ is the unit normal times the sum of the principal curvatures.   In general, $\bH=-\sum_iA(e_i,e_i)$ where $A$ is the second fundamental form and $e_i$ is an orthonormal frame for the tangent space of $\Sigma$; $A(e_i,e_j)= A_{ij} =\nabla^{\perp}_{e_i}e_j$ where $\nabla$ is the Euclidean derivative and `$\perp$' is the component orthogonal to the submanifold.       When $\Sigma$ is noncompact,  $\Sigma_{s,V}$ is replaced by $\Gamma_{s,V}$ where $\Gamma$ is any compact subset of $\Sigma$ containing the support of $V$.    

The submanifold $\Sigma$ is said to be a {\it{minimal}} if 
\begin{align}
\frac{d}{ds}_{s=0}\Vol(\Sigma_{s,V})=0\text{ for all $V$}\, ,
\end{align}
or, equivalently, by \eqr{e:firstvariation}, if  $\bH$ is identically zero.   Thus $\Sigma$ is minimal if and only if it is a critical point for the volume functional.   Since a critical point is not necessarily a minimum, the term minimal is misleading but time-honored.  It is easy to see that being minimal  is equivalent to all the coordinate functions of $\RR^N$ restricted to the submanifold are harmonic with respect to the Laplacian, $\Delta_{\Sigma}$, on the submanifold.    In higher codimension, the minimal surface equation is a complicated system.   

A computation shows that if $\Sigma$ is minimal, then the second derivative of volume is
\begin{align}
\frac{d^2}{ds^2}_{s=0}\Vol  (\Sigma_{s,V})=-\int_{\Sigma}\langle V, L\,V\rangle\, ,
\end{align}
where $L\,V=\Delta_{\Sigma}\,V+\langle A_{ij},\V\rangle\,A_{ij}$ is the so-called second variational (or Jacobi) operator.   This is an operator on the normal bundle of $\Sigma$ and is the Laplacian plus a zeroth order term.  When the submanifold is a hypersurface, this simplifies and becomes $L\,V=\Delta_{\Sigma}\,V+|A|^2\,V$, where $|A|^2$ is the sum of the squares of the principal curvatures.    It simplifies further if one identifies $V$ with its projection $\phi=\langle V,\nn\rangle$ onto the unit normal $\nn$.   Then $L\,\phi=\Delta_{\Sigma}\,\phi+|A|^2\,\phi$.  

A minimal submanifold is {\it{stable}} if it passes the second derivative test
\begin{align}
\frac{d^2}{ds^2}\Vol (\Sigma_{s,V})\geq 0\text{ for all }V\, .
\end{align}
Obviously, if a minimal  surface is area or volume minimizing among competitors with the same boundary, then it is stable as well.  However,  stability is   much more general than being minimizing.
Stability becomes a question about whether the Jacobi operator $L$ is nonnegative or not.   The operator $L$ is much simpler for hypersurfaces and, in particular, 
  it is easy to see that a minimal graph is stable and, more generally, so are multi-valued graphs.   In higher codimension, the question of stability becomes much more complicated because of the vector valued nature of $L$ and  the curvature of the normal bundle.

A classical theorem of Bernstein from 1916 shows that entire (that is where the domain of definition is all of $\RR^2$) minimal graphs in $\RR^3$ are planes.  Whether this is true in higher dimensions became known as the Bernstein problem.   This problem played an important role in the field for decades and  is closely related to   regularity   for area minimizers.   In 1965 and 1966, De Giorgi and Almgren proved the Bernstein theorem for graphs in $\RR^4$ and $\RR^5$.  In 1968, Simons extended the Bernstein theorem to $\RR^6$, $\RR^7$ and $\RR^8$, which was shown to be sharp the next year by Bombieri, De Giorgi and Giusti.  Simons' influential paper introduced the second variation operator and stability to minimal surface theory.  Stability of hypersurfaces was studied by Schoen-Simon-Yau, who showed that as long as the dimension of the hypersurface is at most six and the volume of balls are up to a constant the same as Euclidean balls of the same radius and dimension, then all stable minimal hypersurfaces are planes.    In $\RR^3$ Fischer-Colbrie and Schoen showed the same but without assuming area bounds.    This was also proved independently by Do Carmo and Peng.  Schoen later showed a local version of this that has had huge influence on the development of minimal surfaces in three dimensions.   See \cite{CM2}--\cite{CM4} and \cite{P} for more about minimal surfaces.

The situation is much more complicated in higher codimension where there is no analog of the Bernstein theorem.
A simple argument of Wirtinger from the 1930s, using Stokes' formula, shows that any complex submanifold of $\CC^N$ is volume minimizing among things with the same boundary and, thus, a stable minimal submanifold.   This gives a plethora of area-minimizing, and thus also stable, minimal submanifolds once the codimension is at least two.  Moreover, these examples can have arbitrarily large areas.  
%For example, 
 % the parametrized complex submanifold $z\to (z,z^m)$ is a stable minimal variety that is topologically a plane for each integer $m$.  It has $\Area (B_r\cap \Sigma)\geq C\,m\,r^2$ for $%r\geq 1$. 
Remarkably, Micallef \cite{M} proved a converse in $\RR^4$.    Namely, he showed  that a stable oriented, parabolic minimal surface   in $\RR^4$ is complex for some orthogonal complex structure.  Being parabolic is a conformal property that  holds, for instance, if the volume of balls grows at most quadratically.

\begin{figure}[htbp]
\centering\includegraphics[totalheight=.4\textheight, width=.75\textwidth]{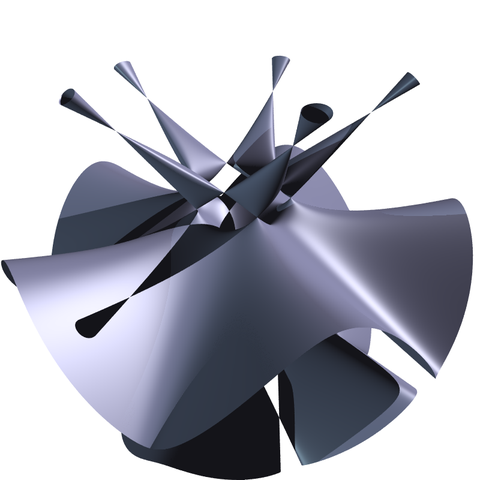}
\caption{A complex curve is a stable minimal surface in the complex plane.}   
  \end{figure}

\section{Motion by mean curvature}

Surface tension is the tendency of fluid surfaces to shrink into the minimum surface area possible.   
Mathematically, the force of surface tension is described by the mean curvature.
A one-parameter family of $n$-dimensional submanifolds $M_t\subset \RR^N$ is said to move by motion by mean curvature if the time derivative of the position vector $x$ 
moves by minus the mean curvature.   That is,
\begin{align}   \label{e:mcfeq1}
 \frac{\partial x}{\partial t} =-\bH\,.
\end{align}
Coordinate functions of the ambient Euclidean space restricted to the evolving submanifolds satisfy the heat equation
\begin{align}   \label{e:mcfeq2}
   \frac{\partial x}{\partial t} =\Delta_{M_t}\,x\, .
\end{align}
This equation is nonlinear since the Laplacian $\Delta_{M_t}$ depends on $M_t$.  Moreover, 
since the submanifolds are evolving, the induced metric   is time-varying so the Laplacian $\Delta_{M_t}$ is also time-varying.     From the first variation formula \eqr{e:firstvariation}, it follows easily that    mean curvature flow   moves in the direction where the volume decreases as fast as possible; thus, mean curvature flow is the negative gradient flow of volume.     The motion is by surface tension.    In higher codimension \eqr{e:mcfeq1} and \eqr{e:mcfeq2} are complicated parabolic systems where much less is known. 

\begin{figure}[htbp]
    \begin{minipage}[t]{0.5\textwidth}
\centering\includegraphics[totalheight=.32\textheight, width=.9\textwidth]{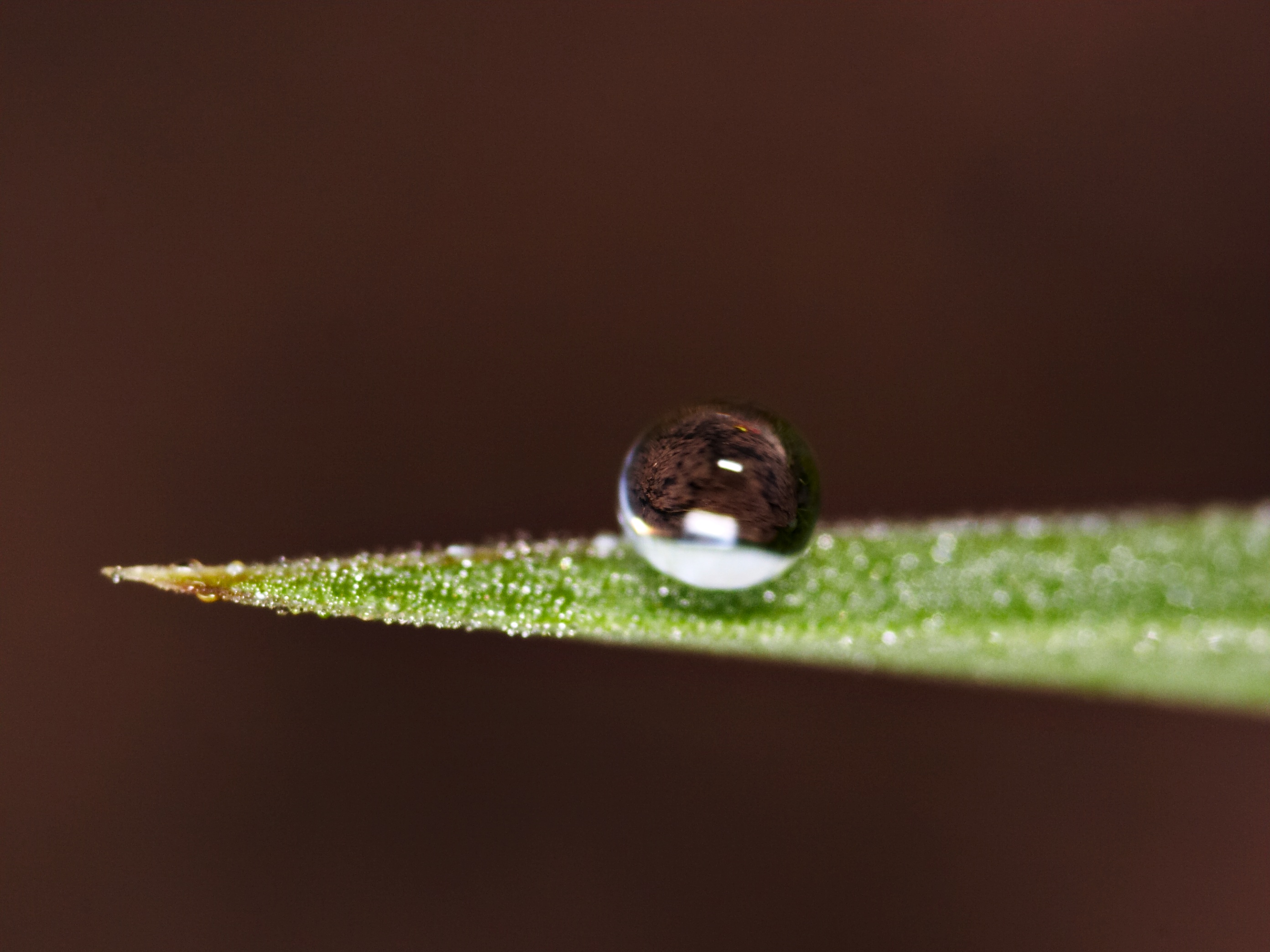}
\caption{Surface tension is the tendency of fluid surfaces to shrink into the minimum surface area possible.   }   
 \end{minipage}\begin{minipage}[t]{0.5\textwidth}
\centering\includegraphics[totalheight=.32\textheight, width=.9\textwidth]{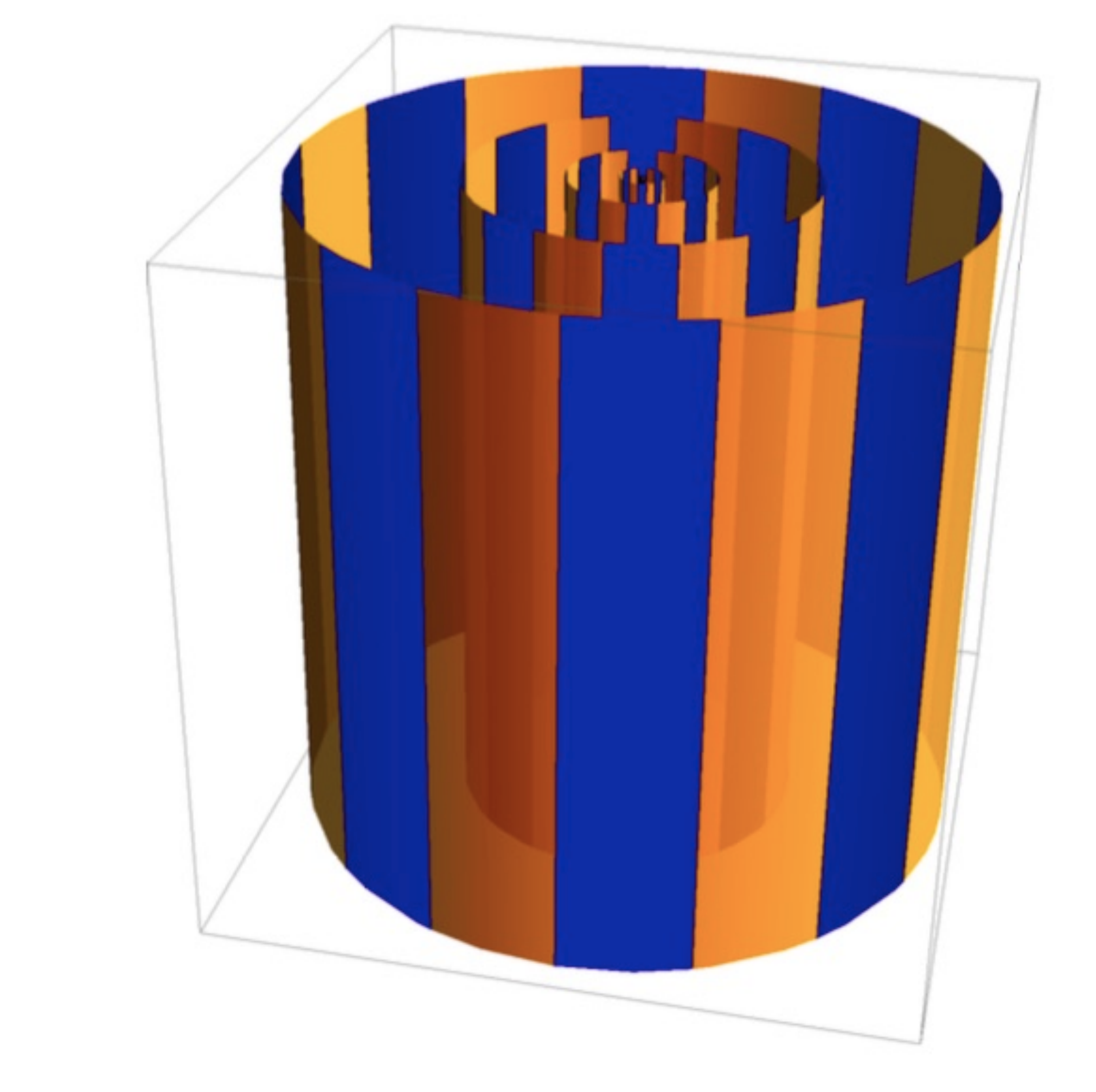}
\caption{Cylinders shrinking homothetically under mean curvature flow.   By \cite{CM5}, shrinking cylinders and spheres are the only stable hypersurface singularities.}   
\end{minipage}
  \end{figure}

Since the coordinate functions on the evolving submanifolds satisfy the heat equation, it follows from the parabolic maximum principle that the evolving submanifolds remain inside the convex hull of the initial submanifold.  A straightforward computation shows that also the function $|x|^2-2\,n\,t$ satisfies the heat equation on the evolving submanifolds.  At the initial time $t=0$ this is nonnegative and therefore, by the parabolic maximum principle, it remains nonnegative as long as the flow exists.  Since we have already seen that $\max_{M_t} |x|^2$  remains bounded under the evolution, it follows that the flow must become extinct in finite time and, thus, singularities occur.   

For a fixed constant $c>0$, rescaling the flow parabolically 
\begin{align}
t\to c\,M_{c^{-2}}=M_{c,t}\, ,
\end{align}
gives a new solution to motion by mean curvature that has the effect that the submanifolds are magnified by the constant $c$.   
If we simultaneously with rescaling also reparametrize time, then we get a rescaled mean curvature flow.  
It is easy to see that such a one-parameter family satisfies the rescaled mean curvature flow equation
\begin{align}
   \frac{\partial x}{\partial t} =\frac{x^{\perp}}{2}-\bH\, .
\end{align}
The rescaled mean curvature flow turns out to be the negative gradient flow for the Gaussian surface area.   

The Gaussian surface area $F$ of an $n$-dimensional submanifold
 $\Sigma^n \subset \RR^N$ is
\begin{align}
	F(\Sigma) = \left( 4\,\pi \right)^{ - \frac{n}{2}} \,  \int_{\Sigma} \e^{ - \frac{|x|^2}{4} } \, .
\end{align}
The constant $\left( 4\,\pi \right)^{ - \frac{n}{2}} $ is a normalization that makes the Gaussian  area one for an $n$-plane through the origin.
Following \cite{CM3}, the entropy $\lambda$ is the supremum of $F$ over all translations and dilations
\begin{align}    
	 \lambda (\Sigma) = \sup_{c,x_0} \, F (c\,\Sigma + x_0)\, .
\end{align}

If $V$ is a normal vector field and $\Sigma_{s,V}$, as before, is the variation $\Sigma_{s,V}=\{x+s\,V(x)\,\,|\,\,x\in\Sigma\}$, then an easy computation shows that
\begin{align}
\frac{d}{ds}_{s=0}F(\Sigma_{s,V})=(4\,\pi)^{-\frac{n}{2}}\int_{\Sigma}\langle V,\bH-\frac{x^{\perp}}{2}\rangle\,\e^{-\frac{|x|^2}{4}}\, .
\end{align}
It follows that the Gaussian surface area $F$ is monotone non-increasing under the rescaled mean curvature flow and constant if and only if
\begin{align}
\bH=\frac{x^{\perp}}{2}\,.
\end{align}
This equation is the shrinker equation and is equivalent to the rescaled flow is static.  Or, also equivalently, the evolution under the mean curvature flow is by rescaling.  That is, a later time slice is exactly like an earlier, just scaled down.    That Gaussian surface area  is monotone under the rescaled flow corresponds  to Huisken's celebrated monotonicity formula, \cite{H}.    From this, it follows also that the entropy is a Lyapunov function for both the mean curvature flow and the rescaled mean curvature flow.   

From Huisken's monotonicity, \cite{H}, and work of Ilmanen, \cite{I}, White, \cite{W}, one knows   that every sequence $c_i\to \infty$ has a subsequence (also denoted by $c_i$) so that $M_{c_i,t}$ converges to a shrinker $M_{\infty,t}$ (so $M_{\infty,t}=\sqrt{-t}\,M_{\infty,-1}$) with $\sup_t\lambda (M_{\infty,t})\leq \sup_t\lambda (M_t)$.  Such a limit is said to be a tangent flow at the  origin.   Similarly one can magnify (blow up) around any other point in space time.   If one doesn't fix  the point around where one blows up, but still looks at limits of a sequence of blow  ups, then the limiting flows are not shrinkers, but even then the limiting flows will exist  for all negative times and are said to be ancient flows.

We have already seen that shrinkers are critical points for the Gaussian area.   When $\Sigma$ is a shrinker, we therefore look at the second derivative.  A calculation (see \cite{CM5}) gives 
\begin{align}
\frac{d^2}{ds^2}_{s=0}F(\Sigma_{s,V})=-(4\,\pi)^{-\frac{n}{2}}\int_{\Sigma}\langle V,L\,V\rangle\,\e^{-\frac{|x|^2}{4}}\, .
\end{align}
Here $L\,V=\cL\,V+\langle A_{ij},V\rangle\,A_{ij}+\frac{1}{2}\,V$ is the Jacobi operator, and $\cL\,V=\Delta_{\Sigma}\,V-\frac{1}{2}\nabla^{\perp}_{x^T}V$ is the Ornstein-Uhlenbeck operator on the normal bundle.   For hypersurfaces, 
there is a similar simplification of the operator $L$ as we saw for the second derivative of volume.  

 For any shrinker, translations and scaling give directions where the Gaussian area decreases, \cite{CM5}, so there are no stable shrinkers in the usual sense.  This corresponds to  $E^{\perp}$ (where $E\in\RR^N$ is a fixed vector) and $\bH=\frac{x^{\perp}}{2}$ being eigenvectors of $L$ with eigenvalues $-\frac{1}{2}$ and $-1$, respectively.  Perturbing by either translation or scaling has the effect of moving the same singularity to a different point in space or time.  However, the singularity is not avoided; it just occurs at another time or place for the flow.   For this reason, we say (\cite{CM5}) that a shrinker is {\it{$F$-stable}} if 
\begin{align}
\frac{d^2}{ds^2}_{s=0}F(\Sigma_{s,V})\geq 0\text{ for all $V$ orthogonal to $\bH$ and all $E^{\perp}$}\, .
\end{align}
Here orthogonal means with respect to the Gaussian inner product on the space of normal vector fields.    For noncompact shrinkers, it turns out that the right notion of stability is that of entropy stability, however, for compact singularities those two notions of stability are the same, \cite{CM5}.  A shrinker is {\it{entropy-stable}} if it is a local minimum for the entropy $\lambda$. Entropy-unstable shrinkers are singularities that can be perturbed away, whereas entropy-stable ones cannot; see \cite{CM5}.  The paper \cite{CIMW} showed that for hypersurfaces round spheres are the shrinkers with smallest entropy.
%\cite{CIMW} conjectured further that round spheres had the least entropy for any closed hypersurface; this was proven by Bernstein-Wang up to dimension $7$ and extended by Zhu to higher dimensions.   
 It is easy to see that spheres and planes are $F$-stable in any codimension.     

\begin{figure}[htbp]
    \begin{minipage}[t]{0.5\textwidth}
\centering\includegraphics[totalheight=.25\textheight, width=.9\textwidth]{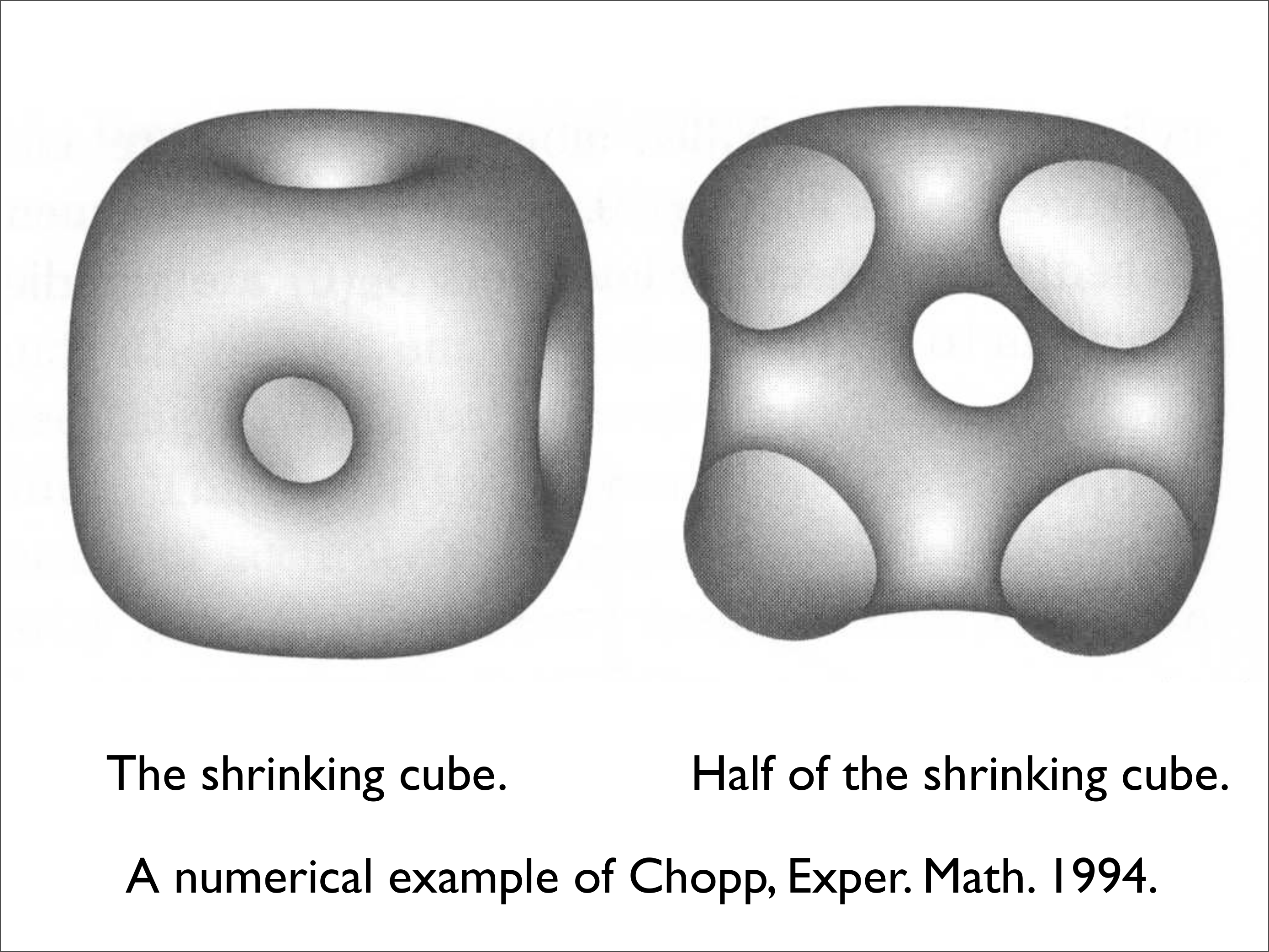} 
\caption{A shrinker in $\RR^3$ found numerically by Chopp.   By \cite{CM5} it is unstable.}   \label{f:chopp1}
 \end{minipage}\begin{minipage}[t]{0.5\textwidth}
\centering\includegraphics[totalheight=.32\textheight, width=.9\textwidth]{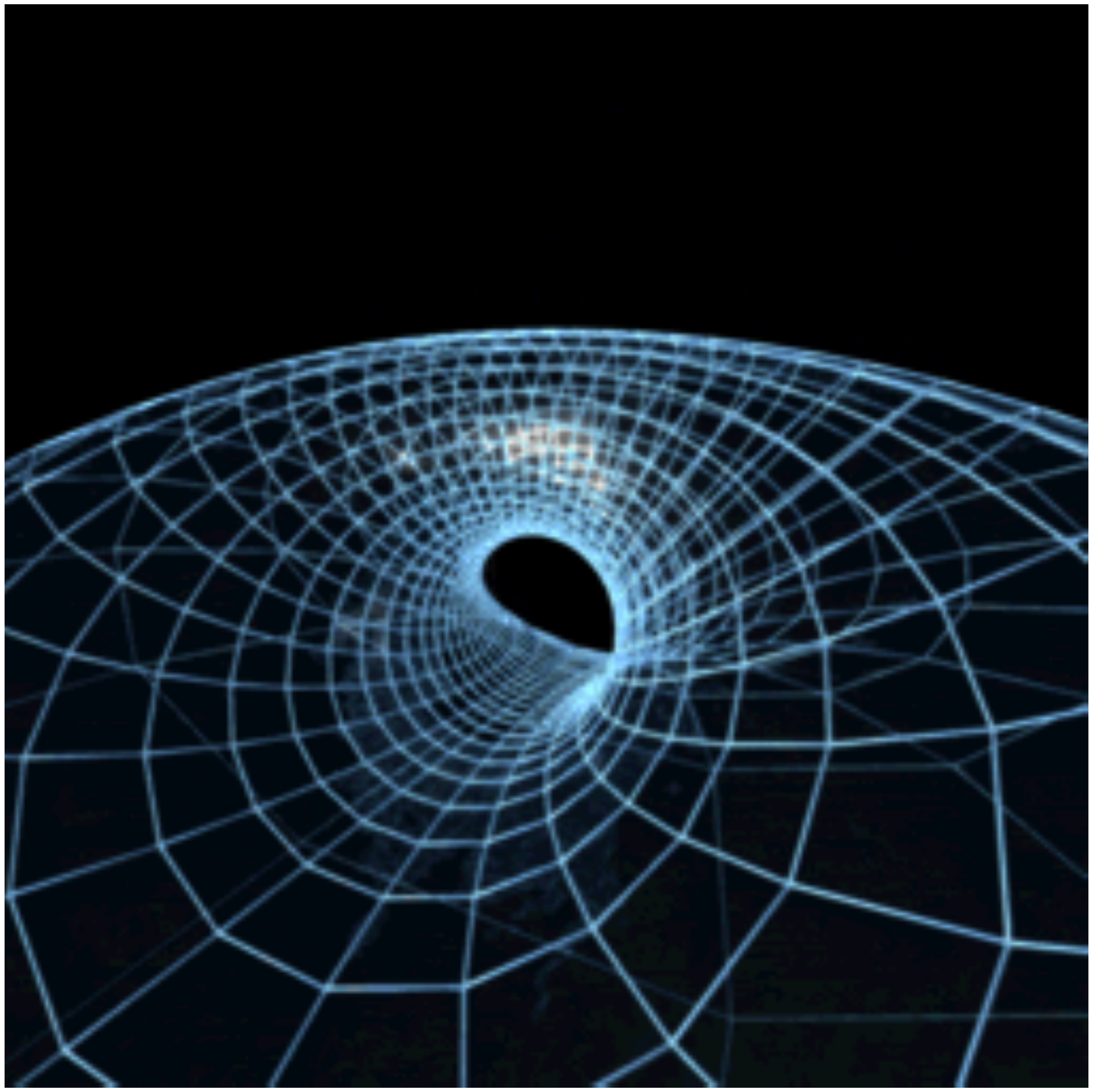}   
\caption{Minimal submanifolds of round spheres are shrinkers.  The illustration is a projection of the Clifford torus (product torus) $\SS^1_{\sqrt{2}}\times \SS^1_{\sqrt{2}}$ that is minimal surface in the round three sphere $\SS^3_2$ of radius $2$ and, thus, a shrinker in $\RR^4$.}   \label{f:cliff}
\end{minipage}
  \end{figure}

Even for hypersurfaces,  examples show that singularities of mean curvature flow are too numerous to classify, see figures \ref{f:chopp1} and \ref{f:shrinker9a}.     The hope is that the generic ones that cannot be perturbed away are much simpler.   
Indeed, in all dimensions, generic singularities (that is, entropy-stable shrinkers) of hypersurfaces moving by mean curvature flow have  been classified in \cite{CM5}.  These are round generalized cylinders $\SS^k_{\sqrt{2\,k}}\times \RR^{n-k}$.   The generic singularities in $\RR^3$ are the sphere $\SS^2_2$, cylinder $\SS^1_{\sqrt{2}} \times \RR$ and plane $\RR^2$. In contrast to the Bernstein theorems for minimal hypersurfaces, this classification of generic singularities holds in every dimension.

 \begin{figure}[htbp]
\centering\includegraphics[totalheight=.4\textheight, width=.75\textwidth]{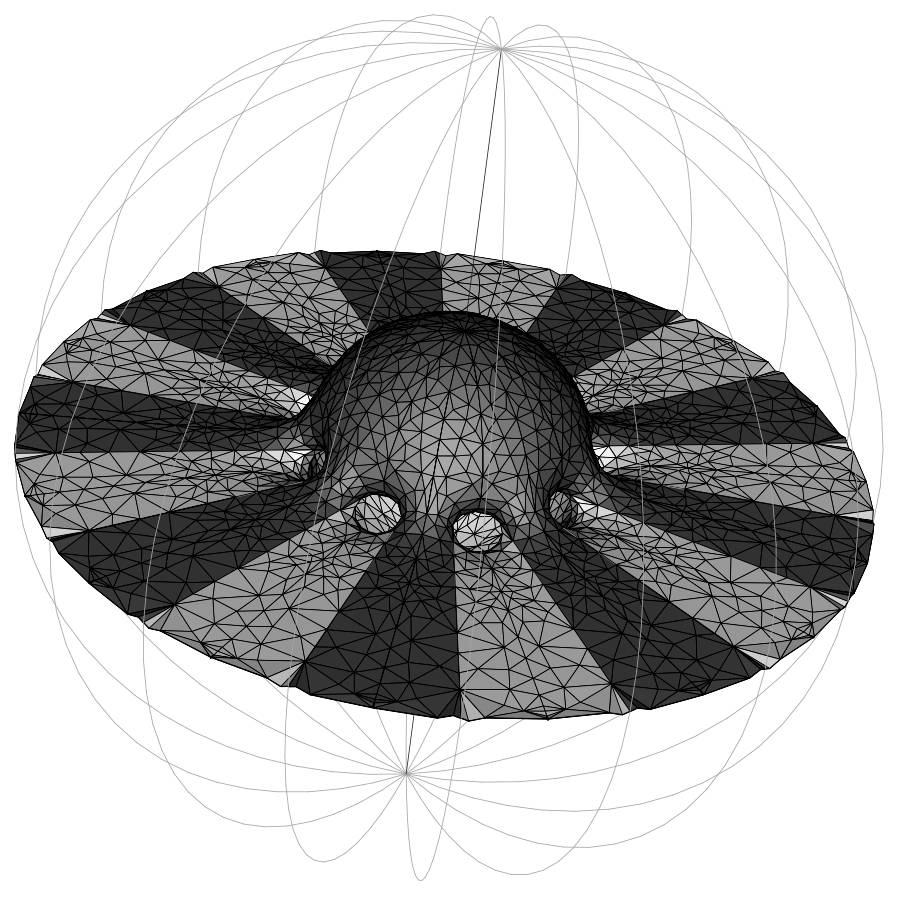}  
\caption{A shrinker in $\RR^3$ shown to exist by Kapouleas, Kleene and M\o ller.  Its existence had been conjectured by Ilmanen using numerics.   By \cite{CM5} it is unstable.}    \label{f:shrinker9a}
  \end{figure}

For the mean curvature flow in higher codimension, we search again for the stable singularities.   Recall that stable singularities are those that are entropy stable, which is equivalent to $F$-stable for closed shrinkers.  When $\Sigma$ is an $F$-stable shrinker diffeomorphic to a sphere, \cite{CM7} shows  that 
\begin{align}	\label{e:becomes}
	\lambda (\Sigma) &<  4 = \e \, \lambda (\SS^2_2) \, .
\end{align}
The sharp constant is unknown, but \eqr{e:becomes} is at most off by a 
  factor of   $\e$.
By \cite{CM7}, similar bounds also hold  for other closed shrinking surfaces where the  entropy bound depends on the genus.     Combined with \cite{CM6}, this implies that any such $F$-stable shrinker, that, a priori, lies in a high dimensional Euclidean space, in fact, lies in a linear subspace of some fixed small dimension.    The sharp bound for the dimension of the linear space is unknown though \cite{CM6} provides sharp dimension bounds in various other important situations.  

There is no analog of \eqr{e:becomes} for minimal surfaces   in $\RR^4$.  Namely, viewing $\RR^4$ as $\CC^2$, one sees that the parametrized complex submanifold $z\to (z,z^m)$ is a stable minimal variety that is topologically a plane for each integer $m$.  It has $\Area (B_r\cap \Sigma)\geq C\,m\,r^2$ for $r\geq 1$.   In contrast, \cite{CM7} implies that  $\Area (B_r\cap\Sigma)\leq C\,(1+\gamma)\,r^2$ for a closed stable $2$-dimensional shrinker  $\Sigma$  of genus $\gamma$.     Similarly, there is no analog of the codimension bound for minimal surfaces.   Indeed, for each $m$, the parametrized surface $z\to (z,z^2,z^3,\cdots,z^{m+1})$ is a stable minimal variety that is topologically a plane.  Its real codimension is $2\,m$ and it is not contained in a proper subspace.   

Once one has the entropy bound in \eqr{e:becomes}, to conclude that stable singularities have low codimension, one needs a result about the number of linearly independent coordinate functions.  
The coordinate functions   on a mean curvature flow produce a linear space of caloric functions, i.e., solutions of the heat equation, that grow at most linearly.
The bound on the codimension is a consequence of a much more general result about polynomial growth caloric functions on an ancient mean curvature flow that has a variety of other useful applications.  

Let $M_t^n \subset \RR^N$ be an ancient mean curvature flow of $n$-dimensional submanifolds with entropies 
$\lambda (M_t)\leq \lambda_0<\infty$.  Recall that ancient flows are solutions that exist for all negative times.   The  space  $\cP_d$ of polynomial growth caloric functions consists of   $u(x,t)$ on $\cup_t M_t\times\{t\}$ so that  $(\partial_t - \Delta_{M_t})\, u=0$ and there exists $C$ depending on $u$ with
  \begin{align}
  	|u(x,t)| \leq C\, (1 + |x|^d+|t|^{\frac{d}{2}}) {\text{ for all }} (x,t) \text{ with }x\in M_t , \, t < 0 \, .
  \end{align}
   The simplest example is when the flow consists of a static (constant in time) hyperplane $\RR^n$.  In this case, $\cP_d(\RR^n)$ consists of polynomials in
   $(t, x_1, \dots x_n)$
 known as the  caloric polynomials and, using the special structure in this case,  it is easy to see that $\dim \cP_d (\RR^n) \approx c_n \, d^n$.  
The paper \cite{CM6} showed  sharp   bounds for $\dim \cP_d$ for all $d\geq 1$ for an ancient flow with $\lambda (M_t) \leq \lambda_0$
\begin{align}
 \dim \cP_d \leq C_n \, \lambda_0 \, d^{n}\, .
 \end{align}
 One remarkable consequence when $d=1$ is a bound for the codimension.  Namely, the flow sits inside a linear subspace of dimension at most $\dim \, \cP_1$,  
  since a linear relation for coordinate functions specifies a hyperplane containing the flow.   
  
 When the manifold $M^n$ is fixed (and does not vary in time), there is a natural subspace of $\cP_d$ consisting of the harmonic functions of polynomial growth of degree at most $d$.  The study of harmonic functions of polynomial growth has a long history in geometry and analysis.  In the 1974, Yau conjectured that these spaces were finite dimensional for manifolds with nonnegative Ricci curvature; this was   proven in \cite{CM1}.
  See \cite{CM8}, \cite{LZ} for results about caloric functions of polynomial growth on such a fixed manifold and the survey \cite{CM9} for other results and the history.

\end{document}